\documentclass[11pt]{amsart}
\usepackage{amssymb,amsmath,amsfonts, tikz, graphicx, subfig, microtype, color, cite, comment}
\usepackage{verbatim,wasysym}
\usepackage[left=2.7cm,right=2.7cm,top=3cm,bottom=3cm]{geometry}
\usepackage{microtype}
\usepackage{color,enumitem,graphicx}
\usepackage[colorlinks=true,urlcolor=blue, citecolor=red,linkcolor=blue,
linktocpage,pdfpagelabels, bookmarksnumbered,bookmarksopen]{hyperref}
\usepackage[hyperpageref]{backref}
\usepackage[english]{babel}

\renewcommand{\S}{\mathbb S}

\renewcommand{\H}{{\mathcal H}}

\newcommand{\dist}{{\rm dist}}
\newcommand{\eps}{\varepsilon}
\newcommand{\half}{\frac{1}{2}}

\newcommand{\T}{{\mathscr T}}

\newcommand{\N}{\mathbb N}

\newcommand{\pnorm}[2][]{\if #1'' \left|#2\right|_p \else \left|#2\right|_{#1} \fi}
\newcommand{\R}{\mathbb R}

\newenvironment{properties}[1]{\begin{enumerate}

}{\end{enumerate}}
\DeclareMathOperator{\supp}{supp}

\newtheorem{lemma}{Lemma}[section]

\newtheorem{proposition}[lemma]{Proposition}
\newtheorem{theorem}[lemma]{Theorem}

\theoremstyle{definition}

\newtheorem{remark}[lemma]{Remark}

\numberwithin{equation}{section}

\title[Quasi-linear nonlocal problems]{Recent progresses in the theory \\ of nonlinear nonlocal problems}

\author[S.\ Mosconi]{Sunra Mosconi}
\author[M.\ Squassina]{Marco Squassina}

\address[S.\ Mosconi]{Dipartimento di Informatica e Matematica
	\newline\indent
	Universit\`a degli Studi di Catania
	\newline\indent
	Viale Andrea Doria 6, 95125 Catania, Italy}
\email{mosconi@dmi.unict.it}

\address[M.\ Squassina]{Dipartimento di Matematica e Fisica 
	\newline\indent
	Universit\`a Cattolica del Sacro Cuore 
	\newline\indent
	Via dei Musei 41, 25121 Brescia, Italy}
\email{marco.squassina@unicatt.it}


\thanks{The authors are members of the Gruppo Nazionale per l'Analisi Matematica, la Probabilit\`a e le loro Applicazioni
	(GNAMPA) of the Istituto Nazionale di Alta Matematica (INdAM)}

\thanks{Bruno Pini Mathematical Analysis Seminar -
	Dipartimento di Matematica, Universit\`{a} di Bologna}

\begin{document}

\begin{abstract}
We overview some recent existence and regularity results in the theory of nonlocal nonlinear problems driven by the fractional $p$-Laplacian.
\end{abstract}

\maketitle


\section{Introduction}

Given $s\in(0, 1)$ and $p>1$, the Gagliardo seminorm of a measurable function $u:\R^N\to \R$ is
\[
[u]=\left(\int_{\R^N\times\R^N}\frac{|u(x)-u(y)|^p}{|x-y|^{N+ps}}\, dx\, dy\right)^{\frac 1 p}.
\]
For a bounded open set $\Omega$ in $\R^N$, the $s$-fractional $p$-Laplacian operator $(- \Delta_p)^s$ (with Dirichlet boundary conditions) is defined as the differential of the functional $u\mapsto \frac{1}{p}[u]^p$,
defined on
\[
W^{s,p}_0(\Omega):=\big\{u\in L^p(\R^N): [u]<\infty, \,\, u=0 \text{ in $\R^N\setminus \Omega$}\big\},
\]
which is a Banach space with respect to the norm $[\cdot]$. Under suitable smoothness conditions on $u$ the operator can be written as
\begin{equation*}
(- \Delta_p)^s\, u(x) = 2 \lim_{\varepsilon \searrow 0} \int_{B_\varepsilon^c(x)} \frac{|u(x) - u(y)|^{p-2}\, (u(x) - u(y))}{|x - y|^{N+sp}}\, dy, \quad x \in \R^N.
\end{equation*}
For several motivations concerning the introduction of this class of nonlinear operators, we refer the reader e.g.\ to \cite{Caffa}.
 The aim of this paper is to review some recent {\em existence} and {\em regularity} results for the weak solutions of the quasi-linear problem
\begin{equation}
\label{ilprob}
	\begin{cases}
		(- \Delta_p)^s u =f(x,u) & \text{in $\Omega$}  \\
		u  = 0  & \text{on $\partial\Omega$},
	\end{cases}
\end{equation}
where $f:\Omega\times\R\to\R$ is a Carath\'eodory function satisfying suitable growth conditions. 
A weak solution $u\in W^{s,p}_0(\Omega)$ of \eqref{ilprob}
satisfies, for every $\varphi\in W^{s,p}_0(\Omega)$, 
\[
\int_{\R^N\times\R^N}\frac{|u(x)-u(y)|^{p-2}(u(x)-u(y))(\varphi(x)-\varphi(y))}{|x-y|^{N+ps}} \,dx\,dy=\int_\Omega f(x,u)\varphi(x) \,dx.
\]
In Section~\ref{regul} we overview the regularity results starting from the linear case
(Section~\ref{linear}), namely the fractional Laplacian and then moving forward
to the nonlinear case (Section~\ref{nonlinearc}), that is problems involving the fractional $p$-Laplacian or even more general operators. In Section~\ref{entire} we consider problems defined on the whole $\R^N$, first discussing the optimal decay rate of optimizers of the fractional Sobolev inequality (Section~\ref{critical-c}) and then
the validity of a Poh\v ozaev identity which allows to derive useful nonexistence results (Section~\ref{poho}). As far as the existence of solutions is concerned, in Section~\ref{exist} we state some results on the variational spectrum
of $(-\Delta_p)^s$ (Section~\ref{eigen}), some existence results for subcritical nonlinearities obtained via Morse theory (Section~\ref{subcritical}) and, finally,
we discuss a nonlinear nonlocal version of the celebrated
Brezis-Nirenberg problem (Section~\ref{Brezisnir}).
\vskip3pt
\noindent
Part of the content of this review paper was  presented by the second author during a {\em Bruno Pini Mathematical Analysis Seminar} 
which was held at the Department of Mathematics of the University of Bologna on October 15, 2015.

\section{Regularity}
\label{regul}
Most of the regularity results needed to deal with quasi-linear problems involving the fractional $p$-Laplacian can be derived from those regarding the model problem
\begin{equation}
\label{nonhom}
\begin{cases}
(-\Delta_p)^s u=f&\text{in $\Omega$},\\
u\equiv 0&\text{in $\R^N\setminus\Omega$},
\end{cases}
\end{equation}
for some $f$ belonging to a suitable function space, e.g. $f\in L^q(\Omega)$ or $f\in C^\alpha(\overline{\Omega})$.
\subsection{The linear case}
\label{linear}
Naturally, the first results on the regularity of nonlocal problems regard the case $p=2$ in \eqref{nonhom}, i.e., the fractional Laplacian. In this case the interior regularity problem is classical, and fine boundary regularity results have been found in \cite{ROS, Grubb1}. These have later been refined in \cite{ROS2, ROS3, Grubb2} and applied in a wide number of semilinear problems involving the fractional Laplacian, see \cite{ROSp, IMS0, IanLiuPerSqu}. 

Moving beyond the fractional Laplacian case, there are essentially two classes of linear nonlocal operators which can be considered. Notice that equation \eqref{nonhom} for $p=2$ can be written in two ways:
\begin{enumerate}
\item
The pointwise formulation
\[
-\int_{\R^N} \frac{u(x+h)+u(x-h)-2u(x)}{|h|^{N+2s}}\, dh= f(x),\qquad x\in \Omega
\]
which is well defined (avoiding the principal value) for regular $u$.
\item
The weak formulation
\[
\int_{\R^N\times\R^N}\frac{(u(x)-u(y))(\varphi(x)-\varphi(y))}{|x-y|^{N+2s}}\, dx\, dy=\int_{\R^N} f\varphi\, dx,\qquad \forall \varphi\in C^\infty_c(\Omega)
\]
which is well defined for $u$ belonging to suitable energy spaces.
\end{enumerate}
The first one gives a meaning to $(-\Delta)^su$ for continuous $u$ through the machinery of viscosity solutions. The second one is variational, since it naturally arises when seeking for solutions through nonlinear analysis techniques. These two formulations are equivalent in most cases, at least for the fractional Laplacian. However, even at the linear level the two notions generalize into very different classes of equations, analogous to divergence versus non-divergence form equations. For example, one can look for viscosity solutions of 
\begin{equation}
\label{ndiv}
\int_{\R^N} K(x, x+h)\bigl(u(x+h)+u(x-h)-2u(x)\bigr)\, dh= f(x),
\end{equation}
or for weak solutions of
\begin{equation}
\label{div}
\int_{\R^N\times\R^N}H(x, y)(u(x)-u(y))(\varphi(x)-\varphi(y))\, dx\, dy=\int_{\R^N} f(x)\varphi(x)\, dx,
\end{equation}
both for suitably singular kernels.
The latter linear nonlocal equation is of type \eqref{ndiv} if and only if
\begin{equation}
\label{nldiv}
H(x, x+h)=H(x-h, x), \quad \text{for all $x, h\in \R^N$}
\end{equation}
holds, which is the nonlocal analogue of the PDE condition ${\rm div} A=0$, ensuring that
\[
{\rm div}(A Du)=A\cdot D^2u.
\]
Condition \eqref{nldiv} certainly holds when $H(x, y)=H(x-y)$, but otherwise \eqref{ndiv} and \eqref{div} (and the respective regularity theory) differ substantially, exactly as in the corresponding local equations. 

\subsection{The nonlinear case}
\label{nonlinearc}
The regularity theory for \eqref{ndiv}, with $K$ in suitable kernel classes is well developed, and starting from \cite{CafSil1, CafSil2}, it has been extended to various corresponding classes of fully nonlinear elliptic integro-differential equations. For related bibliography and most recent results, see \cite{Serra} for interior regularity  and \cite{ROS2} for boundary regularity.

On the contrary, much less is known for equations "in divergence form" as in \eqref{div}, and in most cases different techniques are needed. Since the case $p=2$ of \eqref{nonhom} is of both types, we will describe the most recent regularity result on the case $p\neq 2$ of \eqref{nonhom}, which instead cannot be reduced to a fully nonlinear equation built through problems of the type \eqref{ndiv}. Some of the results we will list are in fact proved for the more general operators $(-\Delta_{K, p})^s$ defined on $W^{s,p}_0(\Omega)$ as
\[
\langle (-\Delta_{K, p})^su, \varphi\rangle:=\int_{\R^N\times\R^N}J_p((u(x)-u(y))(\varphi(x)-\varphi(y))K(x, y)\, dx\, dy, 
\]
with $K$ satisfying 
\[
 \frac{\lambda}{|x-y|^{N+ps}}\le K(x, y)\le \frac{\Lambda}{|x-y|^{N+ps}}\qquad \text{for some $0<\lambda\le \Lambda$}.
\]
For simplicity, we will focus on the case $K\equiv 1$, restating the theorems in this case.

The first regularity result has been proved for $f\equiv 0$ and no "boundary condition" except a controlled behavior at infinity ensuring that the operator $(-\Delta_p)^s$ is well defined. Let $\Omega\subseteq \R^N$ be bounded. We set
\[
\widetilde{W}^{s,p}(\Omega):=\Big\{u\in L^p_{\rm loc}(\R^N):\exists\, U\Supset\Omega \ \text{such that\ }\|u\|_{W^{s,p}(U)}+\int_{\R^N}\frac{|u(x)|^{p-1}}{(1+|x|)^{N+ps}}\, dx<\infty\Big\}.
\]
Then $(-\Delta_p)^s u$ is a well defined continuous linear functional on $W^{s,p}_0(\Omega)$.
\begin{theorem}[Interior regularity, \cite{DKP1}]
Let $u\in \widetilde{W}^{s,p}(\Omega)$ satisfy $(-\Delta_p)^s u=0$ in $\Omega$. Then, $u$ is locally H\"older continuous in $\Omega$.
\end{theorem}

Unfortunately, due to the nonlinear nature of the operator it is not immediate to reduce non-homogeneous equations with vanishing Dirichlet conditions to homogeneous equation with nonvanishing Dirichlet conditions.
However, for nonhomogeneous equations, general fine continuity results when $f$ is a measure are contained in \cite{Kuusi}, whose proofs were modified in \cite{BP} to cater for optimal summability conditions on $f$ (namely $f\in L^q$ for $q>\frac{N}{sp}$) ensuring H\"older regularity.
For boundedly nonhomogeneous equations, the interior regularity proof has been simplified in \cite{IaMosSq}, where continuity up to the boundary for problem \eqref{nonhom} has also been established.

\begin{theorem}[Boundary regularity, \cite{IaMosSq}]
Let $\Omega$ be a bounded domain with $C^{1,1}$ boundary, and $u$ solve \eqref{nonhom} for some $f\in L^\infty(\Omega)$. Then  $u\in C^\alpha(\R^N)$ for some $\alpha>0$, with the estimate $\|u\|_{C^\alpha(\R^N)}\le C(N, p, s, \Omega)\|f\|_{\infty}^{1/(p-1)}$.
\end{theorem}

Recently, in \cite{KKP}, H\"older continuity up to the boundary has also been established for solutions to the fractional obstacle problem arising from the energy
\[
\int_{\R^N\times\R^N}K(x, y)|u(x)-u(y)|^p\, dx\, dy,\qquad \frac{\lambda}{|x-y|^{N+ps}}\le K(x, y)\le \frac{\Lambda}{|x-y|^{N+ps}},
\]
for some $0<\lambda\le \Lambda$. In the framework discussed here (i.e. $\lambda=\Lambda$ and no obstacle), this amounts to remove the $C^{1,1}$ regularity hypothesis on $\partial\Omega$, at least when $f=0$ and $u\in\widetilde{W}^{s,p}(\Omega)\cap C^\alpha(U\setminus\Omega)$ for some $U\Supset \Omega$. Indeed, the more natural measure-theoretic condition
\begin{equation}
\label{mt}
\text{$\exists\  r_0>0$ such that}\qquad  \inf_{x_0\in \partial\Omega}\inf_{0<r<r_0}\frac{|(\R^N\setminus \Omega)\cap B_r(x_0)|}{r^N}>0
\end{equation}
is employed therein. 
Finally, higher Sobolev regularity has also been established in the case $p\ge 2$. We state here one of the results contained in \cite{BL}.

\begin{theorem}[Sobolev regularity, \cite{BL}]
\label{BL}
Suppose $p\ge 2$ and let $u$ solve \eqref{nonhom} for $f\in W^{s,\frac{p}{p-1}}(\R^N)$. Then, for any $0<\theta<1$ it holds
\[
s\le\frac{p-1}{p+1}\quad \Rightarrow \quad u\in W^{\theta s\frac{p+1}{p-1}, p}_{\rm loc}(\Omega),
\]
\[
s>\frac{p-1}{p+1}\quad \Rightarrow\quad u\in W^{\theta s\frac{p+2}{p}, p}_{\rm loc}(\Omega).
\]
\end{theorem}

Let us mention some open problems regarding the model equation \eqref{nonhom}. 
\begin{enumerate}
\item
{\em (Higher H\"older interior regularity).} 
It is not known what is the optimal interior regularity for $u$ solving \eqref{nonhom} under H\"older regularity conditions on $f$. Even in the homogeneous case, i.e.
\begin{equation}
\label{hom}
(-\Delta_p)^s u=0\quad \text{in $\Omega$}, \qquad u\in C^\infty(\R^N\setminus\Omega)\cap L^\infty(\R^N)\subseteq \widetilde{W}^{s,p}(\Omega)
\end{equation}
only unspecified H\"older continuity is known.  Homogeneous examples suggest that one cannot expect more that $C^{s\frac{p}{p-1}}_{\rm loc}(\Omega)$ interior regularity and a natural conjecture would be that any solution to \eqref{hom} satisfies $u\in C^{s(1+\eps)}_{\rm loc}(\Omega)$ for some $ \eps>0$. The same conjecture holds for nonhomogeneous equations with H\"older continuous right hand side.
\item
{\em (Higher Sobolev interior regularity).} Clearly the case $p\in \ ]1, 2[$ is missing in Theorem \ref{BL}. However, comparing to the local case $s=1$, one expects very different results in the singular case. On the other hand, for $p\ge 2$ the assumptions on $f$ in Theorem  \ref{BL} are not optimal (as already stated by the authors in \cite{BL}). For $p\neq 2$, the general picture seems yet to be understood, also because even in the local case $s=1$ optimal differentiability assumptions have only very recently been found in \cite{BS}.
\item
{\em (Boundary regularity)}
The techniques developed in \cite{KKP} probably lead to H\"older regularity up to the boundary under condition \eqref{mt} also for boundedly nonhomogeneous equations like \eqref{nonhom}, but a proof is missing. Notice that $\alpha$-H\"older regularity up to the boundary implies, for the solutions to \eqref{nonhom}, the upper bound
\begin{equation}
\label{bdist}
|u(x)|\leq C(N, p, s, \Omega, f){\rm dist}(x, \R^N\setminus \Omega)^\alpha.
\end{equation}
The latter estimate for $\alpha=s$ (which is the optimal H\"older exponent due to explicit examples) was proved in \cite{IaMosSq} and turned out to be a successful tool to obtain existence for the fractional Brezis-Nirenberg problem, see \cite[Theorem 1.3, case ({\em iv})]{PYMS}. It is then relevant to know what is the minimal condition on $\partial\Omega$ ensuring \eqref{bdist} for $\alpha=s$. The examples in \cite{ROS4} suggest that $\partial\Omega$ being $C^{1,\beta}$, $\beta>0$ can suffice, while for $C^1$ boundaries the optimal exponent in \eqref{bdist} may be strictly less than $s$.

\item
{\em (Higher boundary regularity)}
In \cite{ROS}, one of the results states that for $p=2$ and $C^{1,1}$ domains, any solution to \eqref{nonhom} with $f\in L^\infty$ satisfies $u/{\rm d}^s\in C^\alpha(\overline{\Omega})$ for some $\alpha>0$, where ${\rm d}(x)={\rm dist}(x, \R^N\setminus\Omega)$. This higher boundary regularity is the nonlocal counterpart of $u\in C^{1,\alpha}(\overline{\Omega})$ for solutions of $-\Delta_pu=f\in L^\infty$. It has later been generalized and refined to a suitable class of fully nonlinear operators in \cite{ROS2, ROS3}. Coupled with a Hopf-type lemma, it has important applications to nonlinear analysis, exploited for example in \cite{IanLiuPerSqu, IMS0, PYMS}. No such result is known for \eqref{nonhom}, even if the previously mentioned validity of \eqref{bdist} for $\alpha=s$ suggests the validity of the following estimate in $C^{1,1}$ domains:
\begin{equation}
\label{RC}
\|\frac{u}{{\rm d}^s}\|_{C^\alpha(\overline{\Omega})}\le C(N, p, s, \Omega)\|f\|_\infty^{\frac{1}{p-1}},\qquad \text{ for some $\alpha>0$.}  
\end{equation}
\end{enumerate}

\section{Entire problems}
\label{entire}

Another fundamental tool to deal with problem \eqref{ilprob} is the study of entire solutions (i.e.,  $\Omega=\R^N$) for autonomous nonlinearities, namely
\begin{equation}
\label{ent}
(-\Delta_p)^s u= f(u),\qquad \text{in $\R^N$}. 
\end{equation}
The entire problem indeed arises in many blow-up arguments, and the important case $f(u)=\lambda |u|^{q-2}u$, $q>1$ is a good model problem. In particular we will discuss the {\em critical case}
\begin{equation}
\label{at}
(-\Delta_p)^s u= |u|^{p^*-2}u\qquad \text{in $\R^N$,}\qquad \quad p^*=\frac{Np}{N-ps},
\end{equation} 
since it is the basic tool to understand the failure of the compactness of $W^{s,p}(\R^N)\hookrightarrow L^{p^*}(\R^N)$ when $N>ps$ and related equations. Notice that the class of solutions to \eqref{at} is translation invariant and have the following scale invariance too 
\[
u\mapsto u_t(x)=t^{\frac{N-ps}{p}}u(t x), \qquad t>0.
\]
In the local case $s=1$ the equation itself has a meaning for a wide variety of functions $u$, namely $u\in W^{1,p}_{\rm loc}(\R^N)$ and one may consider solutions with infinite kinetic energy, finding out fine classification results based on stability and/or Morse index, see \cite{Farina} for the model case $f(u)=|u|^qu$. However, due to the non-local nature of the fractional $p$-Laplacian, a global {\em a-priori} growth control on the solution has to be imposed regardless, in order to give a meaning to the equation. The more natural one is
\[
\int_{\R^N}\frac{|u|^{p-1}}{(1+|x|)^{N+ps}}\, dx<+\infty,
\]
but this can be weakened through other notions of solution, see \cite{Kuusi}, or for solutions of \eqref{ent} in other unbounded domains, e.g. $\Omega=\R^N\setminus\{0\}$, see \cite{BMS}. In order to simplify the exposition we will only consider the so called {\em finite  energy solutions}, i.e. those having finite Gagliardo seminorm. To this end, we will frequently state the results in the space
\[
D^{s,q}(\R^N)=\big\{u\in L^{q^*}(\R^N):\int_{\R^N\times\R^N}\frac{|u(x)-u(y)|^q}{|x-y|^{N+qs}}\, dx\, dy<+\infty\big\}.
\]
\subsection{Critical problems}
\label{critical-c}
We now focus on \eqref{at}, always supposing in the following that $N>ps$. Among finite energy solutions one can single out {\em variational solutions}, i.e. those minimizing the Rayleigh quotient
\begin{equation}
\label{Smin}
\S:=\inf_{u\in L^{p^*}(\R^N)\setminus\{0\}}\frac{[u]^p}{\|u\|_{p^*}^p}, \qquad s\in \ ]0, 1[, \quad p>1.
\end{equation}
whose Euler equation is indeed of the form \eqref{at}.
For $s<1$, the existence of minimizers for \eqref{Smin} requires quite a detailed analysis of non-local interactions and has been settle down in \cite{mosqu, MM} using the concentration-compactness method. 

For $s=1$ the variational solutions have been classified in \cite{talenti, aubin} and turned out to be of all the translations, rescalings and constant multiples of 
\[
U(x)=(1+|x|^{p'})^{\frac{p-N}{p}}, \qquad p'=\frac{p}{p-1},
\]
which are since then called the {\em Aubin-Talenti functions}.
It has been a long standing conjecture that any {\em constant sign} solution to \eqref{at} is an Aubin-Talenti function, solved for $p\le 2$ in \cite{Vetois, ClassDino} and only recently for $p>2$ in \cite{sciunziCl}. If no sign condition is imposed on the energy solutions to \eqref{at}, then for $p=2$ there are actually infinitely many, conformally non-equivalent, of them (see \cite{Ding}). 

For $s<1$ the conjectured form of the Aubin-Talenti functions is the following (up to translations, rescaling and constant multiples)
\begin{equation}
\label{Talenti}
U(x)=(1+|x|^{p'})^{\frac{ps-N}{p}},
\end{equation}
at least for minimizers of \eqref{Smin}, however it not even known if these functions actually solve \eqref{at} for $p\neq 2$.   For  $p=2$ most of the previous results have been proved in \cite{Lieb1, Lieb2, CLO} through various methods, all of which employ stereographic projection and/or some conformal symmetry of the problem. None of these approaches seem work for $p\neq 2$, neither will do the ODE methods in \cite{aubin, talenti}, since even after one-dimensional reduction through rearrangement inequalities, the solution solves a fractional integral equation rather than a differential one. Finally, the optimal transport approach of \cite{cordero} may need very deep modifications since it relies on the Brenier formula which is a strongly local result.

A first result enforcing the conjectured form of the Aubin-Talenti functions for general $p\neq 2$ is proved in \cite{BMS}, where the exact asymptotic behavior is determined for minimizers of \eqref{Smin}.

\begin{theorem}[Asymptotic behavior, \cite{BMS}]
\label{sharp-Dec}
Let $N>ps$ and $U$ be any minimizer for \eqref{Smin}. 
Then $U\in L^\infty(\mathbb{R}^N)$, is of constant sign, radially symmetric and monotone around some point, and it holds
\begin{equation*}
\lim_{|x|\to\infty}|x|^{\frac{N-ps}{p-1}}U(x)=U_\infty\in \R\setminus\{0\}.
\end{equation*}
\end{theorem}

A similar asymptotic behavior has later been extended in \cite{MM} to variational solutions of the more general Hardy-Sobolev inequality, and it is therein proved an energy estimate which is still in accordance with the conjectured form \eqref{Talenti}. 

\begin{theorem}[Summability estimate, \cite{MM}]
\label{energydecay}
Let $N>ps$ and $U$ be any minimizer for \eqref{Smin}. Then $u\in D^{s, q}(\R^N)$, $\forall q>\frac{N(p-1)}{N-s}$.
\end{theorem}

As a matter of fact, Theorem \ref{sharp-Dec} turned out to be enough for most nonlinear analysis applications to \eqref{ilprob}. In particular, through a new technique described in \cite{PYMS},  it allows the construction of a suitable family of compactly supported almost minimizers for \eqref{Smin}, see Lemma \ref{fgt} below.

\subsection{Poh\v ozaev identity}
\label{poho}
Let $X\subseteq L^1_{\rm loc}(\R^N)$ be a Banach space which is dilation invariant, i.e. $u\in X\Rightarrow u_\lambda\in X$ for all $\lambda>0$, where $u_\lambda(x):=u(\lambda x)$. A Poh\v ozaev identity for a functional $J\in C^1(X)$ at a critical point $u$ can be stated as follows:
\[
dJ(u)=0\quad \Rightarrow\quad \left.\frac{d}{d\lambda}J(u_\lambda)\right|_{\lambda=1}=0.
\]
Often, due to suitable scaling properties of $J$, the derivative on the right can be easily computed, giving rise to the commonly known Poh\v ozaev identity for $J$ at $u$.

At this level of generality, there is little hope to obtain a Poh\v ozaev identity even assuming smoothness of $\lambda\mapsto J(u_\lambda)$,
the issue being that the curve $\lambda\mapsto \gamma(\lambda):=u_\lambda\in X$ may be very irregular, preventing the application of the chain rule 
\[
 \left.\frac{d}{d\lambda}J(u_\lambda)\right|_{\lambda=1}=\langle dJ(u), \gamma'(1)\rangle.
\]
A classical example \cite{W} of Whitney giving rise to a variety of Morse-Sard type theorems, shows the existence of $J\in C^k(\R^N)$ ($k<N$) and $\gamma\in C^\alpha([0, 1]; \R^N)$ ($\alpha<1$) such that $dJ(\gamma(t))=0$ for all $t\in [0,1]$ while $t\mapsto J(\gamma(t))$ is non constant. The main reason for such a pathological example is tied to the low regularity condition $k<N$, $\alpha<1$ in the previous example, so that in infinite dimensional spaces the situation gets worse, and even requiring $J\in C^\infty$ doesn't suffice to avoid the previous pathology, see e.g. \cite{Kupka}.
A more refined analysis shows that these Whitney type counterexamples can be constructed when some relation between the regularity of $J$ and that of $\gamma$ are satisfied, see \cite{Norton1, Norton2} for some of these kind of results. 

Clearly, the regularity of the curve $\lambda \mapsto \gamma(\lambda)$ is tied to the regularity of the critical point $u$ since, at least formally, $\gamma'(1)=\nabla u\cdot x$, which may not belong to $X$ (either for regularity or for summability issues). Therefore Poh\v ozaev identities are interesting (and useful) to obtain exactly in those cases where the critical points of $J$ may have low regularity, e.g. for semilinear problems involving the $p$-Laplacian.

Historically, many Poh\v ozaev identities have been stated {\em assuming} higher (not to be expected) regularity of the critical point $u$, see e.g. \cite{PS}. Later on, these assumptions have been removed in \cite{GV, DGMS}, obtaining Poh\v ozaev identities without regularity assumptions on $u$ for very general integral functionals of the calculus of variations.

The situation in the fractional case is much less understood, mainly because the optimal boundary and Sobolev regularity of the critical points of the relevant functional is much lower. For the fractional Laplacian $p=2$ in bounded domains however, the picture is quite complete.

\begin{theorem}[Poh\v ozaev identity, \cite{ROSp}]
Let $f\in {\rm Lip}_{\rm loc}(\R)$, with $F(t)=\int_0^tf(\tau)\, d\tau$ and let $\Omega$ a $C^{1,1}$ bounded domain, with exterior normal $n$ and distance function ${\rm d}(x)={\rm dist}(x, \R^N\setminus\Omega)$. Then any bounded solution $u$ to 
\begin{equation}
\label{ros}
\begin{cases}
(-\Delta)^s u=f(u)&\text{in $\Omega$}\\
u\equiv 0&\text{in $\R^N\setminus \Omega$}
\end{cases}
\end{equation}
satisfies the Poh\v ozaev identity
\[
(2s-N)\int_{\Omega}uf(u)\, dx+2N\int_\Omega F(u)=\Gamma(1+s)^2\int_{\partial\Omega}\left|\frac{u}{{\rm d}^s}\right|^2 x\cdot n\, d{\mathcal H}^{N-1}.
\]
where $\Gamma$ is the Gamma function.
\end{theorem}

Notice that $u/{\rm d}^s$ is well defined and continuous on the whole $\overline{\Omega}$ due to the results in \cite{ROS}. For $p\neq 2$ this boundary regularity is still unknown, but a very general argument still provides the validity of at least a Poh\v ozaev {\em inequality}.

\begin{theorem}[Poh\v ozaev inequality, \cite{ROSne}]
Let $f\in {\rm Lip}_{\rm loc}(\R)$, with $F(t)=\int_0^tf(\tau)\, d\tau$ and let $\Omega$ be a bounded star-shaped domain. Then any bounded solution to \eqref{ros} which in addition lies in $W^{1, r}(\Omega)$, satisifes
\[
N\int_\Omega F(u)\, dx\ge \frac{N-ps}{p}\int_{\Omega} uf(u)\, dx.
\]
\end{theorem}

Here we see that higher Sobolev regularity has to be assumed {\em a priori}, which actually is a quite strong hypothesis for small $s$ and $p\neq 2$.
In the entire case $\Omega=\R^N$ the issue of the lack of informations on the boundary regularity is substituted by the problem of the decay of the solution. As far as we know, even in the the case $p=2$, an entire Poh\v ozaev identity with no additional assumptions on the decay of the solution is still missing, but some very close results are available (see \cite{wangzq}). 

Finally, let us remark that one of the main use of Poh\v ozaev identities lies in ruling out {\em sign changing} solutions to suitable nonlinear problems. In the entire case, variants of the moving plane method are well developed for $p=2$, often giving optimal nonexistence results for {\em constant sign} solutions, see \cite{CLO2}.

\section{Existence}
\label{exist}
We now describe some existence results for problems of the from
\begin{equation*}
	\begin{cases}
		(- \Delta_p)^s u =f(x,u) & \text{in $\Omega$},  \\
		u  = 0  & \text{on $\partial\Omega$}.
	\end{cases}
\end{equation*}
In the following we will always assume that  $f:\Omega\times\R\to\R$ is a Carath\'eodory mapping, and set $F(x,t)=\int_0^t f(x,\tau){\rm d}\tau$ for all $(x,t)\in\Omega\times\R$.

\subsection{Eigenvalue problems}
\label{eigen}
The first investigations in the literature dealt with the case
$$
f(x,t)=\lambda |t|^{p-2}t,\quad t\in\R,\,\,\,  \lambda>0,
$$
namely, with the spectral analysis of the operator $(-\Delta_p)^s$. As in the local case, there are many, in principle different, ways to define the spectrum of $(-\Delta_p)^s$. We will focus on the two most used notions, namely the classical spectrum and the Fadell-Rabinowitz one. More precisely, the classical spectrum will be denoted by $\sigma(\Omega)$ and defined as the set of real numbers such that $(-\Delta_p)^s u=\lambda u$ for some nontrivial $u\in W^{s,p}_0(\Omega)$. For the definition of the Fadell-Rabinowitz spectrum $\sigma_{FR}$ we refer to \cite{IanLiuPerSqu}. 
The main properties of these spectra are enlisted below, where the dependencies of the spectra on $N, p, s$ and $\Omega$ is omitted when no confusion can arise.
\[
\sigma_{FR}\subseteq \sigma, \qquad \min\sigma=\min\sigma_{FR}=:\lambda_1>0
\]
\[
\inf\sigma\setminus\{\lambda_1\}=\inf\sigma_{FR}\setminus\{\lambda_1\}=:\lambda_2.
\] 
Some results naturally parallel the existing local nonlinear spectral theory for the $p$-Laplacian operator,  even if in some cases the proofs have to be heavily modified. In some instances, moreover, the theory differs substantially, giving rise to new phenomena caused by the non-local nature of the operator.

\begin{theorem}[The first eigenvalue, \cite{LL}]
The first eigenvalue of $(- \Delta_p)^s$ is positive and isolated, i.e. $\lambda_2>\lambda_1$. Morever its eigenspace is one dimensional and spanned by a nonnegative function.
\end{theorem}

Another result proved in \cite{LL}, which is in contrast to the local case $s=1$, is a strict inequality which holds due the nonlocal nature of the operator, namely that if $u$ is a continuous sign-changing eigenfunction, then
\[
\lambda>\max\big\{\lambda_1(\{u>0\}), \lambda_1(\{u<0\})\big\}, \qquad \forall \lambda\in \sigma(\Omega).
\]

\begin{theorem}[Hong-Krahn-Szego inequality, \cite{BP}]
It holds
\begin{equation}
\label{HKS}
\lambda_2(\Omega)>\left(2\frac{|B|}{|\Omega|}\right)^{\frac{sp}{N}}\lambda_1(B),
\end{equation}
where $B$ is any ball. Moreover, inequality is never attained,  but it is sharp for a sequence of disjoint balls of identical volume arbitrarily far apart.
\end{theorem}
Notice that in the local case $s=1$, equality is attained in \eqref{HKS} for $\Omega$ being the union of two disjoint balls of equal measure.
\begin{theorem}[Weyl law, \cite{meyan}]
Let $N<sp$ and $\Omega$ be a bounded domain. The non-decreasing sequence of eigenvalues $\{\lambda_k\}_{k\ge 1}=\sigma_{FR}$ (counted with multiplicity) satisfies
$$
C_1 k^\frac{sp-N}{N}\leq \lambda_k\leq C_2 k^\frac{Np-N+sp}{N},
$$
for some positive constants $C_1,C_2$ depending on $\Omega$ and on $s,p,N$.
\end{theorem}

\noindent
Finally, let us mention \cite{PSY2}, where the Fucik spectrum of  $(- \Delta_p)^s$ was investigated, namely
the existence of continuous curves of points $(a,b)$ in $\R^2$ such that problem \eqref{ilprob} with 
$$
f(x,t)=b (t^+)^{p-1}-a (t^-)^{p-1},\quad t\in\R,
$$
admits a nontrivial solution in $W^{s,p}_0(\Omega)$.

\subsection{Subcritical problems}
\label{subcritical}
Consider now the problem \eqref{ilprob} with
\begin{equation*}
f(x,t)=\lambda|t|^{p-2}t+g(x,t) 
\end{equation*}
where $\lambda\in\R$ is a parameter and the hypotheses on the reaction term $g$ are the following:
\begin{itemize}
	\item[${\bf H}_1$] $g:\Omega\times\R\to\R$ is a Carath\'eodory mapping, $G(x,t)=\int_0^t g(x,\tau){\rm d}\tau$, and
	\begin{itemize}
		\item[$(i)$] $|g(x,t)|\le a(1+|t|^{r-1})$ a.e. in $\Omega$ and for all $t\in\R$ ($p<r<p^*_s$);
		\item[$(ii)$] $0<\mu G(x,t)\le g(x,t)t$ a.e. in $\Omega$ and for all $|t|\ge R$ ($\mu>p$, $R>0$);
		\item[$(iii)$] $\displaystyle\lim_{t\to 0}\frac{g(x,t)}{|t|^{p-1}}=0$ uniformly a.e. in $\Omega$.
	\end{itemize}
\end{itemize}
Hypotheses ${\bf H}_1$ classify problem \eqref{ilprob} as $p$-{\em superlinear}. Besides, by ${\bf H}_1(iii)$ we have $g(x,0)=0$ a.e. in $\Omega$, so \eqref{ilprob} admits the zero solution for all $\lambda\in\R$. By means of Morse theory and the spectral properties of $(- \Delta_p)^s$, one can prove the existence of a non-zero solution for all $\lambda\in\R$, requiring when necessary additional sign conditions on $G(x,\cdot)$ near zero.

\begin{theorem}[Existence via Morse theory, \cite{IanLiuPerSqu}]
	\label{4main}
	If ${\bf H}_1$ and one of the following hold:
	\begin{itemize}
		\item[$(i)$] $\lambda\notin\sigma_{FR}$;
		\item[$(ii)$] $\lambda\in\sigma_{FR}$ and $G$ is of constant sign in $\Omega\times [-\delta, \delta]$ for some $\delta>0$;
		
	\end{itemize}
	then problem \eqref{ilprob} admits a non-zero solution.
\end{theorem}

\vskip4pt
\noindent
Let us now consider a new set of assumptions for problem \eqref{ilprob}, in order to get multiple solutions.
\begin{itemize}
	\item[${\bf H}_2$]
	\begin{itemize}
		\item[$(i)$] $f(x,t)t\ge 0$ a.e. in $\Omega$ and for all $t\in\R$;
		\item[$(ii)$] $\displaystyle\lim_{t\to 0}\frac{f(x,t)-b|t|^{q-2}t}{|t|^{p-2}t}=0$ uniformly a.e. in $\Omega$ ($b>0$, $1<q<p$);
		\item[$(iii)$] $\displaystyle\limsup_{|t|\to\infty}\frac{f(x,t)}{|t|^{p-2}t}<\lambda_1$ uniformly a.e. in $\Omega$.
	\end{itemize}
\end{itemize}

In the local case $s=1$, the techniques employed to obtain multiplicity results under assumptions ${\bf H}_2$ are connected with fine boundary regularity theory for nonhomogeneous equations. For $s\in \ ]0, 1[$ this theory is well developed  for linear ($p=2$) problems, see the discussion at point (4) of Section 2. Therefore, for the fractional Laplacian hypotheses ${\bf H}_2$ imply multiplicity results through  a constrained minimization argument, adapted to the fractional setting in \cite{IMS0}.  The general scheme of proof for any $p>1$ has been outlined in \cite{IanLiuPerSqu}, where the expected boundary regularity had instead to be assumed. More precisely, the following holds.

\begin{theorem}[Existence of three solutions, \cite{IanLiuPerSqu}]
	\label{5main}
	If hypotheses ${\bf H}_2$ and $\eqref{RC}$ hold, then problem \eqref{ilprob} admits at least three non-zero solutions.
\end{theorem}

\vskip4pt
\noindent
Existence can be obtained also in the asymptotically $p$-linear case, under the following set of hypotheses:
\begin{itemize}
	\item[${\bf H}_3$] 
	\begin{itemize}
			\item[$(i)$] $\displaystyle\lim_{|t|\to\infty}\frac{f(x,t)}{|t|^{p-2}t}=\lambda>0$ uniformly a.e. in $\Omega$;
		\item[$(ii)$] $\displaystyle\lim_{t\to 0}\frac{f(x,t)-b|t|^{q-2}t}{|t|^{p-2}t}=0$ uniformly a.e. in $\Omega$ ($b>0$, $1<q<p$).
	\end{itemize}
\end{itemize}
Clearly, ${\bf H}_3 (ii)$ implies that $f(x,0)=0$ a.e. in $\Omega$, so \eqref{ilprob} admits the zero solution. 

\begin{theorem}[Asymptotically $p$-linear case I, \cite{IanLiuPerSqu}]
	\label{6ex1}
	If ${\bf H}_3$ holds with $\lambda\notin\sigma_{FR}(\Omega)$, then problem \eqref{ilprob} admits at least a non-zero solution.
\end{theorem}

\vskip3pt
\noindent
As a variant, hypotheses ${\bf H}_3$ can be modified as follows:
\begin{itemize}
	\item[${\bf H'}_3$]
	\begin{itemize}
			\item[$(i)$] $\displaystyle\lim_{|t|\to\infty}\frac{f(x,t)}{|t|^{p-2}t}=\lambda>\lambda_1$ uniformly a.e. in $\Omega$;
		\item[$(ii)$] $\displaystyle\lim_{t\to 0}\frac{f(x,t)}{|t|^{p-2}t}=\mu<\lambda_1$ uniformly a.e. in $\Omega$.
	\end{itemize}
\end{itemize}
Then, we have the following multiplicity result.

\begin{theorem}[Asymptotically $p$-linear case II, \cite{IanLiuPerSqu}]
	\label{6mult}
	If ${\bf H'}_3$ holds, then problem \eqref{ilprob} admits at least two non-zero solutions, one non-negative, the other non-positive.
\end{theorem}

\subsection{Brezis-Nirenberg problem}
\label{Brezisnir}

	The solvability conditions on $\lambda\in\R$ for the local semi-linear problem
	\begin{equation*} 
	\begin{cases}
	- \Delta u  =\lambda u+|u|^{2^*-2}u & \text{in $\Omega$,} \\
	u  = 0 & \text{in $\partial\Omega$,}
	\end{cases}
	\end{equation*}	
	 were studied in the celebrated paper \cite{BrNi} by H.\ Brezis, L.\ Nirenberg in the case $\lambda<\lambda_1$ and later on in \cite{capozzi} for $\lambda\ge \lambda_1$. Then, most of the results were extended to the non linear problem in \cite{azorero, arioli, degio}. The most recent results in \cite{degio} ensure the existence of a nontrivial solution if
	 \begin{enumerate}
	 \item
	 {\em (non-resonant case):}
	 $\lambda\notin \sigma_{FR}$ and $N\ge p^2$.
	 \item
	 {\em (resonant case):}
	 $\lambda\in \sigma_{FR}$ and either 
	 $\frac{N^2}{N+1}>p^2$, or  $\frac{N^3+p^3}{N^2+N}>p^2$ and $\partial\Omega$ is $C^{1,\alpha}$.
	 \end{enumerate}
	 
	On this basis, it is natural to consider the Brezis-Nirenberg problem in the nonlocal case.

\begin{theorem}[Nonlocal Brezis-Nirenberg, \cite{PYMS}]
	\label{mountainpass}
	Suppose that one of the following facts holds:
	\begin{enumerate}
		\item
		{\em (non-resonant case):} $\lambda\notin \sigma_{FR}$ and either $N>p^2s$, or $N=p^2s$ and $0<\lambda<\lambda_1$;
		\item
		{\em (resonant case):} $\lambda\in \sigma_{FR}$ and either  $\frac{N^2}{N+s}>p^2s$, or $\frac{N^3+p^3s^3}{N(N+s)}>p^2s$ and $\partial\Omega$ is $C^{1,1}$-smooth.
	\end{enumerate}
	Then the problem 
		\begin{equation*} 
		\begin{cases}
		(- \Delta)_p^s\, u  =\lambda|u|^{p-2}u+|u|^{p^*-2}u & \text{in $\Omega$} \\
		u  = 0 & \text{in $\R^N \setminus \Omega$,}
		\end{cases}
		\end{equation*}
	admits a nontrivial weak solution. 
\end{theorem}	

The main issue with critical problems is the possible lack of the Palais-Smale condition for the corresponding functional. While this loss of compactness is in some sense intrinsic to the problem, compactness can be recovered at suitable energy levels. Therefore one is led to construct linking structures in such a way that the corresponding min-max energy can be efficiently estimated to be in the range of validity of the Palais-Smale condition. This is usually accomplished constructing suitable truncations of the (rescaled) Aubin-Talenti functions (see section 3.1) and precisely quantify the effect of truncation on the relevant norms. 

Since the explicit expression of the Aubin-Talenti functions in the case $s\in \ ]0, 1[$ is missing, the following lemma is the main tool used to prove the previous theorem, and possibly other critical problems. Compare with a similar result used in \cite{mss} in the case $p=2$.
\begin{lemma}
\label{fgt}
Let $N>ps$.  For any $\delta>\eps>0$ there exists a nonnegative $u_{\eps, \delta}\in W^{s,p}_0(B_\delta)$ such that
\[
[u_{\eps, \delta}]^p\le \S^{\frac{N}{ps}}+C\left(\frac{\eps}{\delta}\right)^{\frac{N-ps}{p-1}}\quad \text{and}\qquad 
\|u_{\eps, \delta}\|_{p^*}^{p^*}\ge \S^{\frac{N}{ps}}-C\left(\frac{\eps}{\delta}\right)^{\frac{N}{p-1}},
\]
for a suitable positive constant $C=C(N, p, s)$.
\end{lemma}

\bigskip

\end{document}